\documentclass[12pt]{article}
\usepackage{amsmath}
\usepackage{amsthm}
\usepackage{amsfonts}
\usepackage{amssymb}
\usepackage{esint}
\usepackage{eucal}
\usepackage{mathrsfs}
\usepackage{graphicx,graphics}
\numberwithin{equation}{section}
\usepackage{graphicx,graphics}
\usepackage{pdfsync}
\usepackage{color}
\def \dis {\displaystyle}

\def \confai {-\kern -.5em\rightharpoonup}
\def \cqfd {\hfill$\Box$}
\def \div{\mbox{\rm div}}

\def \al {\alpha}
\def \be {\beta}

\def \de {\delta}

\def \ep {\varepsilon}

\def \Om {\Omega}

\def \ph {\varphi}

\def \si {\sigma}

\def \NN {\mathbb N}

\def \QQ {\mathbb Q}
\def \RR {\mathbb R}

\def \CC {\mathbb C}
\def \beq {\begin{equation}}
\def \eeq {\end{equation}}
\def \ba {\begin{array}}
\def \ea {\end{array}}
\def \bs {\bigskip}

\def \ecart {\noalign{\medskip}}
\newtheorem{Thm}{Theorem}[section]

\newtheorem{Pro}[Thm]{Proposition}
\newtheorem{Lem}[Thm]{Lemma}

\newtheorem{Adef}[Thm]{Definition}

\newtheorem{Arem}[Thm]{Remark}
\newenvironment{Rem}{\begin{Arem}\rm}{\end{Arem}}
\newtheorem{Aexa}[Thm]{Example}

\newtheorem{Anot}[Thm]{Notation}

\def \refe #1.{\eqref{#1})}
\def \reff #1.{figure~\ref{#1}}
\def \refs #1.{Sec.\,\ref{#1}}
\def \refss #1.{Subsection~\ref{#1}}
\def \refD #1.{Definition~\ref{#1}}
\def \refT #1.{Theorem~\ref{#1}}
\def \refL #1.{Lemma~\ref{#1}}
\def \refC #1.{Corollary~\ref{#1}}
\def \refP #1.{Proposition~\ref{#1}}
\def \refPt #1.{Properties~\ref{#1}}
\def \refR #1.{Remark~\ref{#1}}
\def \refE #1.{Example~\ref{#1}}
\def \refN #1.{Notation~\ref{#1}}
\newcounter{marnote}

\usepackage[colorlinks=true, linkcolor=blue]{hyperref} 

\title{Homogenization of linear transport equations.
\\
A new approach.}
\author{Marc Briane
\\
\normalsize Univ Rennes, INSA Rennes,  CNRS, IRMAR - UMR 6625, F-35000 Rennes, France
\\
\normalsize mbriane@insa-rennes.fr
}
\begin{document}
\maketitle
\begin{abstract}
The paper is devoted to a new approach of the homogenization of linear transport equations induced by a uniformly bounded sequence of vector fields $b_\ep(x)$, the solutions of which $u_\ep(t,x)$ agree at $t=0$ with a bounded sequence of $L^p_{\rm loc}(\RR^N)$ for some $p\in(1,\infty)$. Assuming that the sequence $b_\ep\cdot\nabla w_\ep^1$ is compact in $L^q_{\rm loc}(\RR^N)$ ($q$ conjugate of $p$) for some gradient field $\nabla w_\ep^1$ bounded in $L^N_{\rm loc}(\RR^N)^N$, and that there exists a uniformly bounded sequence $\sigma_\ep>0$ such that $\sigma_\ep\,b_\ep$ is  divergence free if $N\!=\!2$ or is a cross product of $(N\!-\!1)$ bounded gradients in $L^N_{\rm loc}(\RR^N)^N$ if $N\!\geq\!3$, we prove that the sequence $\sigma_\ep\,u_\ep$ converges weakly to a solution to a linear transport equation. It turns out that the compactness of $b_\ep\cdot\nabla w_\ep^1$ is a substitute to the ergodic assumption of the classical two-dimensional periodic case, and allows us to deal with non-periodic vector fields in any dimension. The homogenization result is illustrated by various and general examples.
\end{abstract}
\par\bs\noindent
\section{Introduction}
In this paper we study the homogenization of the sequence of linear transport equations indexed by $\ep>0$,
\beq\label{eque}
\left\{\ba{ll}
\dis {\partial u_\ep\over\partial t}-b_\ep\cdot\nabla_x u_\ep=0 & \mbox{in }(0,T)\times\RR^N,\ N\geq 2
\\ \ecart
u_\ep(0,\cdot)=u_\ep^0 & \mbox{in }\RR^N.
\ea\right.
\eeq
where $T>0$ and $p\in[1,\infty]$ with conjugate exponent $q$. Using the DiPerna-Lions transport theory \cite[Corollary II.1]{DiLi}, if for instance $b_\ep$ is a vector field in $L^\infty(\RR^N)^N\cap W^{1,q}_{\rm loc}(\RR^N)^N$ with bounded divergence and the initial condition $u_\ep^0$ is in $L^p(\RR^N)$, then there exists a unique solution $u_\ep(t,x)$ to equation \eqref{eque} in $L^\infty(0,T;L^p(\RR^N))$.
\par
Tartar \cite{Tar} has showed that the homogenization of first-order hyperbolic equations may lead to nonlocal effective equations with memory effects, and E \cite{E} has also obtained from the homogenization of \eqref{eque} effective higher-order hyperbolic equations. Hence, an interesting problem consists in finding sufficient conditions for which the weak limit of the solution $u_\ep$ to equation \eqref{eque} is still a solution to a first-order transport equation.
This type of homogenization result has first been derived in dimension two by Brenier~\cite{Bre} and by Hou, Xin \cite{HoXi}, assuming that $b_\ep(x)=b(x/\ep)$ where $b$ is a divergence free periodic regular vector field. These works have been extended by E \cite[Sec.\,5]{E} when $b_\ep(x)=b(x,x/\ep)$ with $b(x,y)$ divergence free both in $x$ and $y$, and by Tassa \cite{Tas} when there exists a periodic positive regular function $\sigma$ (which is called an invariant measure for $b$) such that
\beq\label{divb}
\div\,(\sigma b)=0\quad\mbox{in }\RR^2.
\eeq
The main assumption of the periodic framework of \cite{Bre,HoXi,E,Tas} is the ergodicity of the flow associated with $b$ (see, {\em e.g.}, \cite[Lect.\,1]{Sin}, or \cite[Chap.\,II, \S\,5]{ReSi}), namely any periodic invariant function by the flow is constant, or equivalently, for any periodic regular function $v$,
\beq\label{erg}
b\cdot\nabla v=0\;\;\mbox{in }\RR^2\;\Rightarrow\;\nabla v=0\;\;\mbox{in }\RR^2,
\eeq
together with $b\neq 0$ in $\RR^2$.
By virtue of the Kolmogorov theorem (see, {\em e.g.}, \cite[Lect.\,11]{Sin} or \cite[Sec.\,2]{Tas}) in dimension two with $b\neq 0$, condition \eqref{erg} is equivalent to
\[
{\langle b_1\rangle\over\langle b_2\rangle}\notin\QQ.
\]
\par
Here, we present a new approach which holds both in the non-periodic framework and in any dimension with a suitable vector field $b_\ep$. The ergodic assumption \eqref{erg} together with $b\neq 0$ is now replaced by the existence of a sequence $w_\ep^1$ in $C^1(\RR^N)$ and $q\in(1,\infty)$ such that
\beq\label{beDw1e}
0< b_\ep\cdot\nabla w_\ep^1\;\to\;\theta_0>0\quad\mbox{strongly in }L^q_{\rm loc}(\RR^N),
\eeq
which is equivalent in the periodic case to the existence of a periodic gradient $\nabla w$ satisfying
\beq\label{bDw1}
b\cdot\nabla w=1\quad\mbox{in }\RR^N.
\eeq
Moreover, the invariant measure $\sigma$ of the periodic case is replaced by a sequence $\sigma_\ep$ satisfying $0<c^{-1}<\sigma_\ep<c$ for some constant $c>1$, and (see Remark~\ref{rem.rec} for an equivalent expression)
\beq\label{siebei}
\div\,(\sigma_\ep\,b_\ep)=0\;\;\mbox{if }N=2\quad\mbox{and}\quad \sigma_\ep\,b_\ep=\nabla w_\ep^2\times\cdots\times \nabla w_\ep^N\;\;\mbox{if }N\geq 3.
\eeq
The case where $\sigma_\ep\,b_\ep$ is only divergence free in dimension $N\geq 3$ remains open.
In this way the vector field $b_\ep$ is naturally associated with the vector field $W_\ep:=(w_\ep^1,\dots, w_\ep^N)$ which induces a global rectification of the field $b_\ep$ in the direction~$e_1$ (see Remark~\ref{rem.rec}). Then, assuming in addition to \eqref{beDw1e}, \eqref{siebei} that $W_\ep$ is uniformly proper (see condition \eqref{upWe} below) and converges both in $C^0_{\rm loc}(\RR^N)^N$ and weakly in $W^{1,N}_{\rm loc}(\RR^N)^N$, we prove (see Theorem~\ref{thm.hom}) that up to a subsequence $\sigma_\ep\,u_\ep$ converges weakly in $L^\infty(0,T;L^p(\RR^N))$ to a solution $v$ to the transport equation
\beq\label{eqvi}
\left\{\ba{ll}
\dis {\partial v\over\partial t}-\xi_0\cdot\nabla_x\left({v\over\sigma_0}\right)=0 & \mbox{in }(0,T)\times\RR^N
\\ \ecart
v(0,\cdot)=v^0 & \mbox{in }\RR^N,
\ea\right.
\eeq
where $\sigma_0$ is the weak-$\ast$ limit of $\sigma_\ep$ in $L^\infty(\RR^N)$, $\xi_0$ is the weak limit of $\sigma_\ep\,b_\ep$ in $L^{N'}_{\rm loc}(\RR^N)^N$ and $v^0$ the weak limit of $\sigma_\ep\,u_\ep^0$ in $L^p(\RR^N)$. Moreover, if $\sigma_\ep$ converges strongly to $\sigma_0$ in $L^1_{\rm loc}(\RR^N)$ (see Remark~\ref{rem.vsiu}) or $u_\ep^0$ converges strongly to $u^0$ in $L^p_{\rm loc}(\RR^N)$, then up to a subsequence $u_\ep$ converges weakly in $L^\infty(0,T;L^p(\RR^N))$ to a solution $u$ to the transport equation
\beq\label{equ1}
\left\{\ba{ll}
\dis {\partial u\over\partial t}-{\xi_0\over\sigma_0}\cdot\nabla_x u=0 & \mbox{in }(0,T)\times\RR^N
\\ \ecart
u(0,\cdot)=u^0 & \mbox{in }\RR^N.
\ea\right.
\eeq
The convergence of $u_\ep$ also turns out to be strong in $L^\infty(0,T;L^2_{\rm loc}(\RR^N))$ if $u_\ep^0$ converges strongly to $u^0$ in $L^p_{\rm loc}(\RR^N)$ with $p>2$ (see the second part of Theorem~\ref{thm.hom}).
\par
The compactness condition \eqref{beDw1e} is the main assumption of Theorem~\ref{thm.hom}. It is equivalent to the compactness of the product $\sigma_\ep\det(DW_\ep)$ which is connected to the vector field $b_\ep$ by \eqref{siebei}. The examples of Section~\ref{s.exa} show that this condition may be satisfied in quite general situations.
\par
Section~\ref{s.res} is devoted to the statement of the main result and to its proof.
Section~\ref{s.exa} deals by three applications of Theorem~\ref{thm.hom}.
In Section~\ref{ss.exa1} we study the case of a diffeomorphism $W_\ep$ on $\RR^2$ such that $\det(DW_\ep)$ is compact in $L^p_{\rm loc}(\RR^2)$ for some $q\in (1,\infty)$. In Section~\ref{ss.per} we extend the periodic case of \cite{Bre,HoXi,E,Tas} with $b_\ep(x)=b(x/\ep)$ and the periodic case of \cite[Sec.\,4]{Bri} on the asymptotic of the flow associated with $b$, in the light of Theorem~\ref{thm.hom} with a periodically oscillating function $\sigma_\ep(x)=\sigma(x/\ep)$ (see Proposition~\ref{pro.per}). In Section~\ref{ss.flow} we consider the case of a diffeomorphism $W_\ep$ which agrees at a fixed time $t$ to a flow $X_\ep(t,\cdot)$ associated with a suitable vector field $a_\ep$ (see Proposition~\ref{pro.ae}). In this general setting assumption~\eqref{beDw1e} holds simply when $\div\,a_\ep$ is compact in $L^q_{\rm loc}(\RR^N)$ for some $q\in(1,\infty)$.
\subsubsection*{Notations}
\begin{itemize}
\item $\left(e_1,\dots,e_N\right)$ denotes the canonical basis of $\RR^N$.
\item $\cdot$ denotes the scalar product in $\RR^N$ and $|\cdot|$ the associated norm.
\item $I_N$ is the unit matrix of $\RR^{N\times N}$, and $R_\perp$ is the clockwise $90^\circ$ rotation matrix in $\RR^{2\times 2}$.
\item For $M\in\RR^{N\times N}$, $M^T$ denotes the transpose of $M$.
\item $Y_N:=\left[0,1\right)^N$, and $\langle f\rangle$ denotes the average-value of a function $f\in L^1(Y_N)$.
\item For any open set $\Om$ of $\RR^N$ and $k\in\NN\cup\{\infty\}$, $C^k_c(\Om)$, respectively $C^k_b(\Om)$, denotes the space of the $C^k$ functions with compact support in $\Om$, respectively bounded in $\Om$.
\item For $k\in\NN\cup\{\infty\}$ and $p\in[1,\infty]$, $C^k_\sharp(Y_N)$ denotes the space of the $Y_N$-periodic functions in $C_k(\RR^N)$, and $L^p_\sharp(Y_N)$ denotes the space of the $Y_N$-periodic functions in $L^p_{\rm loc}(\RR^N)$ ({\em i.e.} in $L^p(K)$ for any compact set $K$ of $\RR^N$).
\item For $u\in L^1_{\rm loc}(\RR^N)$ and $U=(U_j)_{1\leq j\leq d}\in L^1_{\rm loc}(\RR^N)^N$.
\[
\nabla_x u:=\left(\partial_{x_1},\dots,\partial_{x_N}\right)\quad\mbox{and}\quad DU:=\big[\partial_{x_i}U_j\big]_{1\leq i,j\leq d}.
\]
\item For $\xi^1,\dots,\xi^N$ in $\RR^N$, the cross product $\xi^2\times\cdots\times \xi^N$ is defined by
\beq\label{crpr}
\xi^1\cdot(\xi^2\times\dots\times \xi^N)=\det\,(\xi^1,\xi^2,\dots,\xi^N)\quad\mbox{for }\xi^1\in\RR^N,
\eeq
where $\det$ is the determinant with respect to the canonical basis $(e_1,\dots,e_N)$.
\item $o_\ep$ denotes a term which tends to zero as $\ep\to 0$.
\item $C$ denotes a constant which may vary from line to line.
\end{itemize}
\section{The main result}\label{s.res}
Let $W_\ep=(w_\ep^1,\dots,w_\ep^N)$, $\ep>0$, be a sequence of vector fields in $C^1(\RR^N)^N$
which is {\em uniformly proper}, {\em i.e.} for any compact set $K$ of $\RR^N$ there exists a compact set $K'$ of $\RR^N$ satisfying
\beq\label{upWe}
W_\ep^{-1}(K)\subset K'\quad\mbox{for any small enough }\ep>0,
\eeq
and let $W\in C^1(\RR^N)^N$ be such that
\beq\label{cWe}
W_\ep\to W\;\;\mbox{in }C^0_{\rm loc}(\RR^N)^N\quad\mbox{and}\quad W_\ep\rightharpoonup W\;\;\mbox{in }W^{1,N}_{\rm loc}(\RR^N)^N.
\eeq
Let $b_\ep$ be a vector field in $C^0_b(\RR^N)^N\cap W^{1,q}_{\rm loc}(\RR^N)^N$ with bounded divergence and let $\sigma_\ep$ be a positive function in $C^0(\RR^N)\cap W^{1,q}_{\rm loc}(\RR^N)$ satisfying for some constant $c>1$,
\beq\label{siebe}
c^{-1}\leq\sigma_\ep\leq c\quad\mbox{and}\quad
\sigma_\ep\,b_\ep=\left\{\ba{ll}
R_\perp\nabla w_\ep^2 & \mbox{if }N=2
\\ \ecart
\nabla w_\ep^2\times\cdots\times \nabla w_\ep^N & \mbox{if }N\geq 3,
\ea\right.
\quad \mbox{in }\RR^N.
\eeq
Also assume that for $p\in(1,\infty)$ with conjugate exponent $q$, there exists a positive function $\theta_0$ in $C^0(\RR^N)$ such that
\beq\label{cthe}
\theta_\ep:=b_\ep\cdot\nabla w_\ep^1>0\quad\mbox{in }\RR^N \qquad\mbox{and}\qquad \theta_\ep\to\theta_0>0\quad\mbox{strongly in }L^q_{\rm loc}(\RR^N).
\eeq
Finally, assume:
\begin{itemize}
\item either that there exists a constant $B>0$ such that
\beq\label{beB}
|\div\,b_\ep|\leq B\quad\mbox{a.e. in }\RR^N,
\eeq
\item or the regularity condition
\beq\label{rsiebe}
b_\ep\in C^1_b(\RR^N)^N,\quad \sigma_\ep\in C^1(\RR^N)\quad\mbox{and}\quad u_\ep^0\in C^1(\RR^N).
\eeq
\end{itemize}
\begin{Rem}\label{rem.rec}
The definition \eqref{siebe} of $b_\ep$ can be also written for any dimension $N\geq 2$ as the existence of $(N-1)$ gradients $\nabla w_\ep^2,\dots,\nabla w_\ep^N$ satisfying
\beq\label{siebe2}
\forall\,\xi\in\RR^N,\quad \si_\ep\,b_\ep\cdot\xi=\det\left(\xi,\nabla w_\ep^2,\dots,\nabla w_\ep^N\right).
\eeq
In dimension $N\geq 3$ this is exactly the definition of the cross product $\nabla w_\ep^2\times\cdots\times\nabla w_\ep^N$ (see~\eqref{crpr}). In dimension $N=2$ this means exactly that $\si_\ep\,b_\ep=R_\perp\nabla w_\ep^2$, which is equivalent to
\beq\label{divsiebe}
\div\,(\sigma_\ep\,b_\ep)=0\;\;\mbox{in }\RR^2.
\eeq
However, in dimension $N\geq 3$ condition \eqref{siebe} is stronger than $\sigma_\ep\,b_\ep$ divergence free.
\par
The definition \eqref{siebe} of $b_\ep$ and the definition \eqref{cthe} of $\theta_\ep$ are equivalent to the global rectification of the field $b_\ep$ by the diffeomorphism $W_\ep$
\beq\label{DWebe}
DW_\ep^T\,b_\ep=\theta_\ep\,e_1\quad\mbox{in }\RR^N,
\eeq
in the direction $e_1$ with the compact range $\theta_\ep$.

\end{Rem}
\par
Then, we have the following homogenization result.
\begin{Thm}\label{thm.hom}
Let $T>0$, let $p\in(1,\infty)$ and let $u_\ep^0$ be a bounded sequence in $L^p(\RR^N)$. Assume that conditions \eqref{upWe} to \eqref{cthe} together with \eqref{beB} or \eqref{rsiebe} hold true.
Let $u_\ep$ be the solution to the transport equation \eqref{eque} and set $v_\ep:=\sigma_\ep\,u_\ep$.
Then, up to a subsequence $v_\ep$ converges weakly in $L^\infty(0,T;L^p(\RR^N))$ to a solution $v$ to the transport equation
\beq\label{eqv}
\left\{\ba{ll}
\dis {\partial v\over\partial t}-\xi_0\cdot\nabla_x\left({v\over\sigma_0}\right)=0 & \mbox{in }(0,T)\times\RR^N
\\ \ecart
v(0,\cdot)=v^0 & \mbox{in }\RR^N,
\ea\right.
\eeq
where (${\rm Cof}$ denotes the cofactors matrix)
\beq\label{xi0cofDW}
\xi_0={\rm Cof}\,(DW)\,e_1\in C^0(\RR^N)^N,
\eeq
\beq\label{xisiv0}
\sigma_\ep\,b_\ep\rightharpoonup \xi_0\;\;\mbox{in }L^{N'}_{\rm loc}(\RR^N)^N,\quad
\sigma_\ep\rightharpoonup\sigma_0\;\;\mbox{in }L^\infty(\RR^N)\,\ast,\quad
\sigma_\ep\,u_\ep^0\rightharpoonup v^0\;\;\mbox{in }L^p(\RR^N).
\eeq
Moreover, if in addition $b_\ep\in W^{1,p/(p-2)}_{\rm loc}(\RR^N)^N$ with $p>2$ and the sequence $u_\ep^0$ converges strongly to $u^0\in L^p_{\rm loc}(\RR^N)$ with $\sigma_0\in W^{1,\infty}(\RR)$ and $\xi_0\in L^\infty(\RR^N)^N\cap W^{1,p/(p-2)}_{\rm loc}(\RR^N)^N$, then $u_\ep$ converges strongly in $L^\infty(0,T;L^2_{\rm loc}(\RR^N))$ to the solution $u$ to the transport equation
\beq\label{equ}
\left\{\ba{ll}
\dis {\partial u\over\partial t}-{\xi_0\over\sigma_0}\cdot\nabla_x u=0 & \mbox{in }(0,T)\times\RR^N
\\ \ecart
u(0,\cdot)=u^0 & \mbox{in }\RR^N.
\ea\right.
\eeq
\end{Thm}
\begin{Rem}\label{rem.uni}
If in Theorem~\ref{thm.hom} we assume in addition that $\sigma_0$ is in $W^{1,\infty}(\RR^N)$ and $\xi_0$ belongs to $L^\infty(\RR^N)^N\cap W^{1,q}_{\rm loc}(\RR^N)^N$, then by virtue of \cite[Corollary\,II.1]{DiLi} there exists a unique solution $v$ to the transport equation \eqref{eqv}.
\end{Rem}
\begin{Rem}\label{rem.vsiu}
In addition to the conditions \eqref{upWe} to \eqref{cthe} assume that $\sigma_\ep$ converges strongly in $L^1_{\rm loc}(\RR^N)$ to $\sigma_0\in W^{1,q}_{\rm loc}(\RR^N)$. Then, we have $v=\sigma_0\,u$ and $v^0=\sigma_0\,u^0$ where $u^0$ is the weak limit of $u_\ep^0$ in $L^p(\RR^N)$, which implies that equation \eqref{eqv} is equivalent to equation \eqref{equ}. Therefore, $u_\ep$ converges weakly in $L^\infty(0,T;L^p(\RR^N))$ to a solution $u$ to the transport equation \eqref{equ}.
\end{Rem}
To prove Theorem~\ref{thm.hom} we need the following $L^p$-estimate.
\begin{Lem}\label{lem.Lpeue}
Let $b_\ep\in L^\infty(\RR^N)^N\cap W^{1,q}_{\rm loc}(\RR^N)^N$ with bounded divergence be such that
\begin{itemize}
\item either estimate \eqref{beB} holds true,
\item or both conditions \eqref{siebe} and \eqref{rsiebe} hold true.
\end{itemize}
Then, there exists a constant $C>0$ such that for any $u_\ep^0\in L^p(\RR^N)$ with $p\in[1,\infty)$, the solution $u_\ep$ to equation~\eqref{eque} satisfies the estimate
\beq\label{Lpeue}
\|u_\ep(t,\cdot)\|_{L^p(\RR^N)}\leq C\,\|u_\ep^0\|_{L^p(\RR^N)}\quad\mbox{for a.e. }t\in(0,T),
\eeq
\end{Lem}
\par\noindent
{\bf Proof of Theorem~\ref{thm.hom}.}
First of all, note that by  \eqref{siebe} and \eqref{cthe} we have
\beq\label{detDWe}
\det (DW_\ep)=\sigma_\ep\,\theta_\ep>0\quad\mbox{in }\RR^N.
\eeq
This combined with property \eqref{upWe} and Hadamard-Caccioppoli's theorem \cite{Cac} (or Hadamard-L\'evy's theorem) implies that $W_\ep$ is a $C^1$-diffeomorphism on $\RR^N$. Moreover, since by \eqref{detDWe} $\det(DW_\ep)$ is positive and by \eqref{cWe} $W_\ep$ converges weakly in $W^{1,N}_{\rm loc}(\RR^N)^N$, by virtue of M\"uller's theorem \cite{Mul} $\det(DW_\ep)$ weakly converges to $\det(DW)$ in $L^1_{\rm loc}(\RR^N)$. Hence, passing to the limit in \eqref{detDWe} together with the strong convergence \eqref{cthe} of $\theta_\ep$, the weak convergence \eqref{xisiv0} of $\sigma_\ep$ and the boundedness \eqref{siebe} of $\sigma_\ep$ we get that
\beq\label{detDW}
\det (DW)=\sigma_0\,\theta_0\geq c^{-1}\,\theta_0>0\quad\mbox{a.e. in }\RR^N,
\eeq
which taking into account the continuity of $DW$ and $\theta_0$ implies that $\det (DW)>0$ in $\RR^N$.
Moreover, again by the uniform character of \eqref{upWe} $W$ is a proper mapping.
Therefore,  $W$ is also a $C^1$-diffeomorphism on $\RR^N$.
\par
The weak formulation of equation \eqref{eque} is that for any function $\phi\in C^1_c([0,T)\times\RR^N)$,
\beq\label{weque}
\int_0^T\!\!\int_{\RR^N}u_\ep\,{\partial\phi\over\partial t}\,dx\,dt+\int_{\RR^N}u_\ep^0(x)\,\phi(0,x)\,dx=\int_0^T\!\!\int_{\RR^N}u_\ep\,\div\left(\phi\,b_\ep\right)dx\,dt.
\eeq
Using a density argument with $\sigma_\ep\in W^{1,q}_{\rm loc}(\RR^N)$, we can replace the test function $\phi$ by $\sigma_\ep\,\ph$ for any $\ph\in C^1_c([0,T)\times\RR^N)$. This combined with the divergence free of $\sigma_\ep\,b_\ep$ leads us to the new formulation
\beq\label{eqsieue}
\int_0^T\!\!\int_{\RR^N}\,\sigma_\ep\,u_\ep\,{\partial\ph\over\partial t}\,dx\,dt+\int_{\RR^N}\sigma_\ep(x)\,u_\ep^0(x)\,\ph(0,x)\,dx
=\int_0^T\!\!\int_{\RR^N}u_\ep\,\sigma_\ep\,b_\ep\cdot\nabla_x\ph\,dx\,dt.
\eeq
We pass easily to the limit in the left hand-side of \eqref{eqsieue}. The delicate point comes from the right-hand side of \eqref{eqsieue}.
\par
By the $L^p$-estimate \eqref{Lpeue} of Lemma~\ref{lem.Lpeue} combined with the uniform boundedness of $\sigma_\ep$ in \eqref{siebe} there exists a subsequence, still denoted by $\ep$, such that $v_\ep=\sigma_\ep\,u_\ep$ converges weakly to some function $v$ in $L^\infty(0,T;L^p(\RR^N))$.
\par
Let $\psi\in C^1_c([0,T)\times\RR^N)$ the support of which is contained in some compact set $[t_0,t_1]\times K$ of $[0,T)\times\RR^N$, and define
\beq\label{phe}
\ph_\ep(t,x):=\psi(t, W_\ep(x))\quad\mbox{for }(t,x)\in(0,T)\times\RR^N,
\eeq
so that $\nabla_x\ph_\ep(t,x):=DW_\ep(x)\nabla_y\psi(t,y)$.
Hence, making the change of variables $y=W_\ep(x)$ and using \eqref{DWebe} we deduce that
\beq\label{vebe}
\ba{l}
\dis \int_0^T\!\!\int_{\RR^N}v_\ep(t,x)\,b_\ep(x)\cdot\nabla_x\ph_\ep(t,x)\,dx\,dt=\int_0^T\!\!\int_{W_\ep^{-1}(K)}v_\ep(t,x)\,b_\ep(x)\cdot\nabla_x\ph_\ep(t,x)\,dx\,dt
\\ \ecart
\dis =\int_0^T\!\!\int_Kv_\ep(t,W_\ep^{-1}(y))\,\theta_\ep(W_\ep^{-1}(y))\,e_1\cdot\nabla_y\psi(t,y)\,\det(DW_\ep^{-1})(y)\,dy\,dt.
\ea
\eeq
First, using successively the H\"older inequality combined with the $L^p$-estimate \eqref{Lpeue}, the inclusion \eqref{upWe} and the $L^q$-strong convergence \eqref{cthe} of $\theta_\ep$, we have 
\[
\ba{l}
\dis \left|\,\int_0^T\!\!\int_{K}v_\ep(t,W_\ep^{-1}(y))\,(\theta_\ep-\theta_0)(W_\ep^{-1}(y))\,e_1\cdot\nabla_y\psi(t,y)\,\det(DW_\ep^{-1})(y)\,dy\,dt\,\right|
\\ \ecart
\dis \leq C_\psi\!\!\int_0^T\!\!\left(\int_{K}\!\!\big|v_\ep(t,W_\ep^{-1}(y))\big|^p\det(DW_\ep^{-1})(y)\,dy\right)^{1\over p}\!\!
\left(\int_{K}\!\!\big|(\theta_\ep-\theta_0)(W_\ep^{-1}(y))\big|^q\det(DW_\ep^{-1})(y)\,dy\right)^{1\over p}\!\!dt
\\ \ecart
\dis \leq C_\psi\int_0^T\|v_\ep(t,\cdot)\|_{L^p(K')}\|\theta_\ep-\theta_0\|_{L^q(K')}\,dt=o_\ep,
\ea
\]
which implies that
\[
\ba{l}
\dis \int_0^T\!\!\int_{K}v_\ep(t,W_\ep^{-1}(y))\,\theta_\ep(W_\ep^{-1}(y))\,e_1\cdot\nabla_y\psi(t,y)\,\det(DW_\ep^{-1})(y)\,dy\,dt
\\ \ecart
\dis \int_0^T\!\!\int_{K}v_\ep(t,W_\ep^{-1}(y))\,\theta_0(W_\ep^{-1}(y))\,e_1\cdot\nabla_y\psi(t,y)\,\det(DW_\ep^{-1})(y)\,dy\,dt+o_\ep.
\ea
\]
Next, by the uniform convergence~\eqref{cWe}
\[
\nabla_y\psi(t,W_\ep(x))\to\nabla_y\psi(t,W(x))\quad\mbox{in }C^0_{\rm loc}([0,T]\times\RR^N).
\]
Then, making the inverse change of variables $x=W_\ep^{-1}(y)$ together with \eqref{upWe} and using the weak convergence of $v_\ep$ to $v$ in $L^\infty(0,T;L^p(\RR^N))$, we have
\[
\ba{l}
\dis \int_0^T\!\!\int_{K}v_\ep(t,W_\ep^{-1}(y))\,\theta_0(W_\ep^{-1}(y))\,e_1\cdot\nabla_y\psi(t,y)\,\det(DW_\ep^{-1})(y)\,dy\,dt
\\ \ecart
\dis = \int_0^T\!\!\int_{K'}v_\ep(t,x)\,\theta_0(x)\,e_1\cdot\nabla_y\psi(t,W_\ep(x))\,dx\,dt
=\int_0^T\!\!\int_{K'}v(t,x)\,\theta_0(x)\,e_1\cdot\nabla_y\psi(t,W(x))\,dx\,dt+o_\ep.
\ea
\]
Let $\ph\in C^1_c([0,T)\times\RR^N)$ and define similarly to \eqref{phe}
\[
\ph(t,x):=\psi(t, W(x))\quad\mbox{for }(t,x)\in[0,T)\times\RR^N,
\]
so that $\nabla_x\ph(t,x):=DW(x)\nabla_y\psi(t,y)$.
Therefore, passing to the limit in \eqref{vebe} we obtain that
\beq\label{cvebe}
\ba{l}
\dis \int_0^T\!\!\int_{\RR^N}v_\ep(t,x)\,b_\ep(x)\cdot\nabla_x\ph_\ep(t,x)\,dx\,dt
\\ \ecart
\dis =\int_0^T\!\!\int_{\RR^N}v(t,x)\,\theta_0(x)\,\big(DW(x)^T\big)^{-1}e_1\cdot\nabla_x\ph(t,x)\,dx\,dt+o_\ep.
\ea
\eeq
On the other hand, using \eqref{DWebe}, \eqref{siebe} and the Murat-Tartar div-curl lemma in $L^{N'}$-$L^N$ (see, {\em e.g.}, \cite[Th\'eor\`eme~2]{Mur}) with convergences \eqref{cWe}, \eqref{cthe}, \eqref{xisiv0} we get that
\beq\label{cDWesiebe}
DW_\ep^T(\sigma_\ep\,b_\ep)=\sigma_\ep\,\theta_\ep\,e_1\rightharpoonup DW^T\xi_0=\sigma_0\,\theta_0\,e_1\quad\mbox{weakly in }L^1_{\rm loc}(\RR^N).
\eeq
This combined with \eqref{detDW} yields equality \eqref{xi0cofDW}.
Convergences  \eqref{cvebe} and \eqref{cDWesiebe} imply that
\[
\int_0^T\!\!\int_{\RR^N}v_\ep\,b_\ep\cdot\nabla_x\ph_\ep\,dx\,dt\;\mathop{\longrightarrow}_{\ep\to 0}\;
\int_0^T\!\!\int_{\RR^N}{v\over\sigma_0}\,\xi_0\cdot\nabla_x\ph\,dx\,dt.
\]
Finally, passing to the limit in formula \eqref{eqsieue} with $\ph_\ep$, it follows that for any $\ph\in C^1_c([0,T)\times\RR^N)$,
\[
\int_0^T\!\!\int_{\RR^N}\,v\,{\partial\ph\over\partial t}\,dx\,dt+\int_{\RR^N}v^0(x)\,\ph(0,x)\,dx
=\int_0^T\!\!\int_{\RR^N}{v\over\sigma_0}\,\xi_0\cdot\nabla_x\ph\,dx\,dt,
\]
which taking into account that $\xi_0$ is divergence free yields the weak formulation of the desired limit equation \eqref{eqv}.
This concludes the proof of the first part of Theorem~\ref{thm.hom}.
\par
Now, assume in addition that $b_\ep\in W^{1,p/(p-2)}_{\rm loc}(\RR^N)^N$ with $p>2$ and $u_\ep^0$ converges strongly to $u^0$ in $L^p(\RR^N)$ with $\sigma_0\in W^{1,\infty}(\RR^N)$ and $\xi_0\in L^\infty(\RR^N)^N\cap W^{1,p/(p-2)}_{\rm loc}(\RR^N)^N$.
By \cite[Theorem\,II.3 and Corollary\,II.1]{DiLi} $u_\ep^2$ is the unique solution to the equation \eqref{eque} with initial condition $(u_\ep^0)^2$, or equivalently, for any $\phi\in C^1_c([0,T)\times\RR^N)$,
\[
\int_0^T\!\!\int_{\RR^N}u_\ep^2\,{\partial\phi\over\partial t}\,dx\,dt+\int_{\RR^N}(u_\ep^0)^2(x)\,\phi(0,x)\,dx=\int_0^T\!\!\int_{\RR^N}u_\ep^2\,\div\left(\phi\,b_\ep\right)dx\,dt,
\]
Replacing $u_\ep$ by $u_\ep^2$ in the first part of Theorem~\ref{thm.hom} and using the strong convergence of $u_\ep^0$ we get that the sequence $\sigma_\ep\,u_\ep^2$ converges weakly in $L^\infty(0,T;L^{p/2}(\RR^N))$ to the solution $w$ to the transport equation
\beq\label{eqw}
\left\{\ba{ll}
\dis {\partial w\over\partial t}-\xi_0\cdot\nabla_x\left({w\over\sigma_0}\right)={\partial w\over\partial t}-{\xi_0\over\sigma_0}\cdot\nabla_x w+{\xi_0\cdot\nabla\sigma_0\over\sigma_0^2}\,w=0
& \mbox{in }(0,T)\times\RR^N
\\ \ecart
w(0,\cdot)=\sigma_0\,(u^0)^2 & \mbox{in }\RR^N.
\ea\right.
\eeq
Note that by virtue of \cite[Corollary\,II.1]{DiLi} the solution $w$ to equation \eqref{eqw} is unique due to the regularities $\sigma_0\in W^{1,\infty}(\RR^N)$, $\xi_0\in L^\infty(\RR^N)^N\cap W^{1,p/(p-2)}_{\rm loc}(\RR^N)^N$ with divergence free. Moreover, again by \cite[Theorem II.3 and Corollary\,II.1]{DiLi} $v^2$ is the unique solution to the equation induced by \eqref{eqv}
\[
\left\{\ba{ll}
\dis {\partial(v^2)\over\partial t}-{\xi_0\over\sigma_0}\cdot\nabla_x(v^2)+2\,{\xi_0\cdot\nabla\sigma_0\over\sigma_0}\,v^2=0 & \mbox{in }(0,T)\times\RR^N
\\ \ecart
v^2(0,\cdot)=(\sigma_0\, u^0)^2 & \mbox{in }\RR^N,
\ea\right.
\]
or equivalently, for any $\phi\in C^1_c([0,T)\times\RR^N)$,
\[
\ba{l}
\dis \int_0^T\!\!\int_{\RR^N}v^2\,{\partial\phi\over\partial t}\,dx\,dt+\int_{\RR^N}(\sigma_0\, u^0)^2(x)\,\phi(0,x)\,dx
\\ \ecart
\dis =\int_0^T\!\!\int_{\RR^N}v^2\,\div\left(\phi\,{\xi_0\over\sigma_0}\right)dx\,dt+\int_0^T\!\!\int_{\RR^N}2\,v^2\,{\xi_0\cdot\nabla\sigma_0\over\sigma_0^2}\,\phi\,dx\,dt.
\ea
\]
Replacing the test function $\phi$ by $\ph/\sigma_0$ by a density argument, it follows that for any function $\ph\in C^1_c([0,T)\times\RR^N)$,
\[
\ba{l}
\dis \int_0^T\!\!\int_{\RR^N}{v^2\over\sigma_0}\,{\partial\ph\over\partial t}\,dx\,dt+\int_{\RR^N}\sigma_0(x)\,(u^0)^2(x)\,\ph(0,x)\,dx
\\ \ecart
\dis =\int_0^T\!\!\int_{\RR^N}v^2\,\div\left(\ph\,{\xi_0\over\sigma_0^2}\right)dx\,dt+\int_0^T\!\!\int_{\RR^N}2\,v^2\,{\xi_0\cdot\nabla\sigma_0\over\sigma_0^3}\,\ph\,dx\,dt
\\ \ecart
\dis =\int_0^T\!\!\int_{\RR^N}{v^2\over\sigma_0}\,\div\left(\ph\,{\xi_0\over\sigma_0}\right)dx\,dt+\int_0^T\!\!\int_{\RR^N}{v^2\over\sigma_0}\,{\xi_0\cdot\nabla\sigma_0\over\sigma_0^2}\,\ph\,dx\,dt,
\ea
\]
which shows that $v^2/\sigma_0$ is also a solution to equation \eqref{eqw}. By uniqueness we thus get that $w=v^2/\sigma_0$.
Similarly, the solution $u$ to equation \eqref{equ} agrees with $v/\sigma_0$.
Finally, using these two equalities we have for any compact set $K$ of $\RR^N$,
\[
\ba{l}
\dis \int_0^T\!\!\int_{K}\sigma_\ep(u_\ep-u)^2\,dx\,dt=\int_0^T\!\!\int_{K}(\sigma_\ep\,u_\ep^2-2\,\sigma_\ep\,u_\ep\,u+\sigma_\ep\,u^2)\,dx\,dt
\\ \ecart
\dis \mathop{\longrightarrow}_{\ep\to 0}\;\int_0^T\!\!\int_{K}(w-2\,v\,u+\sigma_0\,u^2)\,dx\,dt=0,
\ea
\]
which concludes the proof of Theorem~\ref{thm.hom}. \cqfd
\par\bs\noindent
{\bf Proof of Lemma~\ref{lem.Lpeue}.}
If the uniform boundedness \eqref{beB} of $\div\,b_\ep$ is satisfied, then using the estimate (17) of \cite[Proposition\,II.1]{DiLi} for the solution to the regularized equation of \eqref{eque} and the lower semi-continuity of the $L^p$-norm ($p<\infty)$ we get estimate \eqref{Lpeue}.
\par
Otherwise, assume that conditions \eqref{siebe} and \eqref{rsiebe} hold true.
Using the regularity of the data the proof is based on an explicit expression of the solution to equation \eqref{eque} from the flow $Y_\ep$ associated with the vector field $b_\ep$ by
\beq\label{Yebe}
\left\{\ba{ll}
\dis {\partial Y_\ep(t,x)\over\partial t}=b_\ep(Y_\ep(t,x)), & t\in\RR
\\ \ecart
Y_\ep(0,x)=x\in\RR^d.
\ea\right.
\eeq
Let $u_\ep^0$ be a function in $C^1(\RR^N)^N\cap L^p(\RR^N)$.
It is classical that the regular solution $u_\ep$ to the transport equation \eqref{eque} is given by
\beq\label{ueXe}
u_\ep(t,x)=u_\ep^0(Y_\ep(t,x))\quad\mbox{for }(t,x)\in[0,T]\times\RR^N.
\eeq
Let $t\in[0,T]$. Making the change of variables combined with the semi-group property of the flow
\[
y=Y_\ep(t,x)\;\Leftrightarrow\;x=Y_\ep(-t,y),
\]
we get that
\beq\label{Lpeue01}
\int_{\RR^N}\big|u_\ep^0(Y_\ep(t,x))\big|^p\,dx=\int_{\RR^N}\big|u_\ep^0(y)\big|^p\,\big|\det(D_y Y_\ep(-t,y))\big|\,dy.
\eeq
Moreover, by \eqref{Yebe} and the Liouville formula we have for any $(\tau,y)\in\RR\times\RR^N$,
\[
\det(D_y Y_\ep(\tau,y))=\exp\left(\int_0^{\tau}(\div\,b_\ep)(Y_\ep(s,y))\,ds\right).
\]
However, since by \eqref{siebe} $\sigma_\ep\,b_\ep$ is divergence free, we have
\[
\ba{l}
\dis \int_0^{\tau}(\div\,b_\ep)(Y_\ep(s,y))\,ds=-\int_0^{\tau}\left({\nabla\sigma_\ep\cdot b_\ep\over\sigma_\ep}\right)(Y_\ep(s,y))\,ds
\\ \ecart
\dis =-\int_0^{\tau}{\partial\over\partial s}\big(\ln \sigma_\ep(Y_\ep(s,y))\big)\,ds=\ln\left({\sigma_\ep(y)\over\sigma_\ep(Y_\ep(\tau,y))}\right).
\ea
\]
This combined with the boundedness of $\sigma_\ep$ in condition \eqref{siebe} implies that
\[
\forall\,(\tau,y)\in\RR\times\RR^N,\quad 0<\det(D_y Y_\ep(\tau,y))={\sigma_\ep(y)\over\sigma_\ep(Y_\ep(\tau,y))}\leq c^2.
\]
Hence, we deduce from \eqref{Lpeue01} that
\[
\int_{\RR^N}|u_\ep(x)|^p\,dx=\int_{\RR^N}\big|u_\ep^0(Y_\ep(t,x))\big|^p\,dx\leq c^2\int_{\RR^N}\big|u_\ep^0(y)\big|^p\,dy,
\]
which yields the desired estimate \eqref{Lpeue}. This concludes the proof of Lemma~\ref{lem.Lpeue}. \cqfd
\section{Examples}\label{s.exa}
The purpose of this section is to illustrate the homogenization of the transport equation \eqref{eque} by various oscillating fields $b_\ep$ which satisfy the assumptions of Theorem~\ref{thm.hom}. It means giving examples of diffeomorphism $W_\ep$ on $\RR^N$ satisfying the rectification \eqref{DWebe} of the vector field $b_\ep$ where the sequence $\theta_\ep>0$ is compact in $L^q_{\rm loc}(\RR^N)$ for some $q\in(1,\infty)$.
\subsection{First example}\label{ss.exa1}
Let $\al_\ep,\al\in C^1(\RR)$ be such that for some constant $c>0$,
\beq\label{ae}
\al_\ep\to\al\;\;\mbox{in }C^0_{\rm loc}(\RR),\quad \al_\ep'\geq c\;\;\mbox{in }\RR,\quad\al_\ep'\to\al'\;\;\mbox{in }L^2_{\rm loc}(\RR),
\eeq
and let $\be_\ep,\be\in C^1(\RR)$ be such that for some constant $C>0$,
\beq\label{bee}
\be_\ep\to\be\;\;\mbox{in }C^0_{\rm loc}(\RR),\quad |\be_\ep|\leq C\;\;\mbox{in }\RR,\quad \be_\ep'\mbox{ is bounded in }L^\infty_{\rm loc}(\RR),
\eeq
Consider the vector field $W_\ep\in C^1(\RR^N)^N$ defined by
\beq\label{Wecj}
W_\ep(x):=\begin{pmatrix}\al_\ep(x_1)\,\exp\big\{\be_\ep(\al_\ep(x_1)\al_\ep(x_2))\big\},\al_\ep(x_2)\,\exp\big\{\!-\!\be_\ep(\al_\ep(x_1)\al_\ep(x_2))\big\}\end{pmatrix},\quad x\in\RR^2,
\eeq
which is based on the characterization of the holomorphic mappings on $\CC^2$ with constant Jacobian \cite{Nis}.
The gradient of $W_\ep$ is given by
\[
\left\{\ba{l}
\dis \nabla w_\ep^1(x)=\exp\big\{\be_\ep(\al_\ep(x_1)\al_\ep(x_2))\big\}\begin{pmatrix}
\al_\ep'(x_1)\big(1+\al_\ep(x_1)\al_\ep(x_2)\be_\ep'(\al_\ep(x_1)\al_\ep(x_2))\big)
\\ \ecart
\al_\ep'(x_2)\al_\ep^2(x_1)\be_\ep'(\al_\ep(x_1)\al_\ep(x_2)\big)
\end{pmatrix}
\\ \ecart
\dis \nabla w_\ep^2(x)=\exp\big\{\!-\!\be_\ep(\al_\ep(x_1)\al_\ep(x_2))\big\}\begin{pmatrix}
-\al_\ep'(x_1)\al_\ep^2(x_2)\be_\ep'(\al_\ep(x_1)\al_\ep(x_2))
\\ \ecart 
\al_\ep'(x_2)\big(1-\al_\ep(x_1)\,\al_\ep(x_2)\be_\ep'(\al_\ep(x_1)\al_\ep(x_2))\big)
\end{pmatrix}.
\ea\right.
\]
Also define $b_\ep:=R_\perp\nabla w_\ep^2$ and $\sigma_\ep:=1$, so that conditions \eqref{siebe} and \eqref{beB} are fulfilled.
\par
By \eqref{ae} and \eqref{bee} we have
\[
\ba{c}
\dis W_\ep(x)\to W(x):=\begin{pmatrix}\al(x_1)\,\exp\big\{\be(\al(x_1)\al(x_2))\big\},\al(x_2)\,\exp\big\{\!-\!\be(\al(x_1)\al(x_2))\big\}\end{pmatrix}
\quad\mbox{in }C^0_{\rm loc}(\RR^2),
\\ \ecart
W_\ep\rightharpoonup W \quad\mbox{in }H^1_{\rm loc}(\RR^2),
\ea
\]
so that conditions \eqref{cWe} is satisfied, and
\beq\label{jacWe}
b_\ep\cdot\nabla w_\ep^1(x)=\det(DW_\ep)(x)=\al_\ep'(x_1)\,\al_\ep'(x_2)\to \al'(x_1)\,\al'(x_2)\quad\mbox{in }L^2_{\rm loc}(\RR^2),
\eeq
so that condition \eqref{cthe} is satisfied with $p=2$.
Moreover, since by \eqref{ae}
\[
\forall\,t\in\RR,\quad |\al_\ep(t)-\al_\ep(0)|\geq c\,|t|,
\]
the sequence $\al_\ep(0)$ converges, and $\be_\ep$ is uniformly bounded in $\RR$, condition \eqref{upWe} holds for $W_\ep$.
\par
Note that the oscillations of the drift $b_\ep$ in equation \eqref{eque} are only due to the oscillations of the sequence $\be'_\ep$ which does not appear in the convergence \eqref{jacWe} of the Jacobian.
\subsection{The periodic case}\label{ss.per}
This section extends the periodic framework of \cite{Bre,HoXi,E,Tas} and \cite[Corollary~4.4]{Bri}.
\par
Let $W=(w^1,\dots,w^N)$ be a vector field in $C^2(\RR^N)^N$, and let $M$ be a matrix in $\RR^{N\times N}$ such that
\beq\label{siW}
\big(x\mapsto W(x)-Mx\big)\in C^1_\sharp(Y_N)^N\quad\mbox{and}\quad \sigma:=\det(DW)>0\;\;\mbox{in }\RR^N.
\eeq
Consider the periodic vector field $b\in C^1_\sharp(Y^N)^N$ defined by
\beq\label{sib}
\sigma\,b:=\left\{\ba{ll}
R_\perp\nabla w^2 & \mbox{if }N=2
\\ \ecart
\nabla w^2\times\cdots\times \nabla w^N & \mbox{if }N\geq 3.
\ea\right.
\eeq
\par
We have the following result.
\begin{Pro}\label{pro.per}
Let $u_\ep^0\in C^1(\RR^N)$ be a bounded sequence in $L^p(\RR^N)$ with $p\in(1,\infty)$.
Assume that conditions \eqref{siW} and \eqref{sib} hold true.
Then, the vector fields $W_\ep$ and $b_\ep$ defined by
\beq\label{Webe}
W_\ep(x):=\ep\,W\left({x\over\ep}\right)\quad\mbox{and}\quad b_\ep(x):=b\left({x\over\ep}\right)\quad\mbox{for }x\in\RR^N,
\eeq
satisfy the assumptions of Theorem~\ref{thm.hom}.
\par\noindent
Moreover, for any sequence $u_\ep^0$ in $L^p(\RR^N)$ such that $\sigma(x/\ep)\,u_\ep^0$ converges weakly to $v^0$ in $L^p(\RR^N)$, the solution $u_\ep$ to equation \eqref{eque} is such that $\sigma(x/\ep)\,u_\ep$ converges weakly in $L^\infty(0,T;L^p(\RR^N))$ to the solution $v$ to the equation \eqref{eqv} with $\sigma_0=\langle\sigma\rangle$ and $\xi_0=\langle\sigma\,b\rangle$.
\end{Pro}
\noindent
{\bf Proof of Proposition~\ref{pro.per}.}
By the quasi-affinity of the determinant (see, {\em e.g.}, \cite[Sec.~4.3.2]{Dac}) and by \eqref{siW} we have
\[
\det(M)=\det\,\langle DW\rangle=\big\langle \det(DW)\big\rangle >0,
\]
and by \eqref{Webe} there exists a constant $C>0$ such that
\beq\label{eWeM}
\forall\,x\in\RR^N,\quad |W_\ep(x)-Mx| \leq C\ep,
\eeq
which implies condition \eqref{upWe}.
Moreover, estimate \eqref{eWeM} and the uniform bounded of $DW_\ep$ imply easily the convergences \eqref{cWe} with the limit $W(x):=Mx$.
\par
On the other hand, the definitions \eqref{siW} of $W$, $\sigma$ and the definition \eqref{sib} of $b$ show clearly that condition \eqref{siebe} and the regularity \eqref{rsiebe} hold true.
Moreover, we have
\[
\theta:=b\cdot\nabla w^1={\det(DW)\over\sigma}=1\quad\mbox{in }\RR^N,
\]
which implies \eqref{cthe} since $\theta_\ep(x)=\theta(x/\ep)=1$.
\par
Moreover, let $u_\ep^0$ be a sequence in $L^p(\RR^N)$ such that $\sigma(x/\ep)\,u_\ep$ converges weakly to $v^0$ in $L^p(\RR^N)$. By virtue of Theorem~\ref{thm.hom} combined with Remark~\ref{rem.uni} and using the weak limit of a periodically oscillating sequence, the sequence $\sigma(x/\ep)\,u_\ep$ converges weakly in $L^p(\RR^N)$ to the solution $v$ to the equation \eqref{eqv} with $\sigma_0=\langle\sigma\rangle$ and $\xi_0=\langle\sigma\,b\rangle$.
The proof of Proposition~\ref{pro.per} is now complete.
\cqfd
\subsection{The dynamic flow case}\label{ss.flow}
In this section we construct a sequence $W_\ep$ from a dynamic flow associated with a suitable but quite general sequence of vector fields $a_\ep$.
\par
Let $a_\ep,a$ be vector fields in $C^1(\RR^N)^N$  such that
\beq\label{cae}
a_\ep \to a\;\;\mbox{in }C^0_{\rm loc}(\RR^N)^N,\quad a_\ep\rightharpoonup a\;\;\mbox{in }W^{1,\infty}_{\rm loc}(\RR^N)^N\ast,
\eeq
and for some constant $A>0$,
\beq\label{aeA}
|a_\ep|+|\div\,a_\ep|\leq A\quad\mbox{in }\RR^N.
\eeq
Also assume that there exists $q\in(1,\infty)$ such that
\beq\label{divae}
\div\,a_\ep\to \div\,a\quad\mbox{in }L^q_{\rm loc}(\RR^N).
\eeq
Consider the dynamic flow $X_\ep$ associated with the vector field $a_\ep$ defined by
\beq\label{Xeae}
\left\{\ba{ll}
\dis {\partial X_\ep(t,x)\over\partial t}=a_\ep(X_\ep(t,x)), & t\in\RR
\\ \ecart
X_\ep(0,x)=x\in\RR^d,
\ea\right.
\eeq
and let $X$ be the limit flow associated with the limit vector field $a$.
\par
Then, from any sequence of flows $X_\ep$ we may derive a general sequence of vector fields $b_\ep$ inducing the homogenization of the transport equation \eqref{eque}.
\begin{Pro}\label{pro.ae}
Let $u_\ep^0$ be a bounded sequence in $L^p(\RR^N)$ with $p\in(1,\infty)$.
Assume that conditions \eqref{cae}, \eqref{aeA}, \eqref{divae} hold true. For a fixed $t>0$, define the vector field $W_\ep:=X_\ep(t,\cdot)$ from $\RR^N$ into $\RR^N$, and the vector field $b_\ep$ by \eqref{siebe} with $\sigma_\ep=1$. Then, the sequences $W_\ep$ and $b_\ep$ satisfy the assumptions of Theorem~\ref{thm.hom}.
\par\noindent
Moreover, for any sequence $u_\ep^0$ converging weakly to $u^0$ in $L^p(\RR^N)$, the solution $u_\ep$ to equation~\eqref{eque} converges weakly in $L^\infty(0,T;L^p(\RR^N))$ to a solution $u$ to the equation \eqref{equ} where $\sigma_0=1$ and $\xi_0={\rm Cof}\,(D_x X(t,x))\,e_1$.
\end{Pro}
\begin{Rem}
There is a strong correspondance  between the conditions \eqref{cae}-\eqref{aeA} and \eqref{divae} satisfied by the vector field $a_\ep$, and respectively the conditions \eqref{cWe} and \eqref{cthe} satisfied by the vector fields $W_\ep$ and $b_\ep$.
\end{Rem}
\noindent
{\bf Proof of Proposition~\ref{pro.ae}.}
First of all, conditions \eqref{siebe} and \eqref{beB} are straightforward, since $\sigma_\ep=1$ and $b_\ep$ is divergence free.
Fix $T>0$. By \eqref{aeA} we have
\beq\label{eXe}
\forall\,t\in[0,T],\ \forall\,x\in\RR^N,\quad |X_\ep(t,x)-x|\leq A\,T,
\eeq
so that the uniform property \eqref{upWe} is satisfied by $W_\ep$.
\par
Let $K$ be a compact set of $\RR^N$. Again by \eqref{eXe} there exists a compact set $K'$ of $\RR^N$ such that
\beq\label{K'}
\big\{X_\ep(t,x),\,(t,x)\in[0,T]\times K\big\}\subset K'.
\eeq 
Let $\de>0$. Since $a_\ep$ converges uniformly to $a$ in $K'$ and $a\in C^1(\RR^N)$ is $k$-Lipschitz in $K'$ for some $k>0$, we have for any small enough $\ep>0$ and for any $t\in[0,T]$, for any $x,y\in K$,
\[
\ba{ll}
\big|X_\ep(t,x)-X_\ep(t,y)\big| & \dis \leq |x-y|+\int_0^t \big|a_\ep(X_\ep(s,x))-a_\ep(X_\ep(s,y))\big|\,ds
\\ \ecart
& \dis \leq \de+|x-y|+k\int_0^t\big|X_\ep(s,x)-X_\ep(s,y)\big|\,ds.
\ea
\]
Hence, by Gronwall's inequality (see, {\em e.g.}, \cite[Sec.\,17.3]{HSD}) we get that for any small enough $\ep>0$,
\[
\forall\,t\in[0,T],\ \forall\,x,y\in K,\quad |X_\ep(t,x)-X_\ep(t,y)\big|\leq(\de+|x-y|)\,e^{kt},
\]
which by \eqref{aeA} implies that for any small enough $\ep>0$,
\[
\forall\,s,t\in[0,T],\ \forall\,x,y\in K,\quad |X_\ep(s,x)-X_\ep(t,y)\big|\leq A\,|s-t|+(\de+|x-y|)\,e^{kt},
\]
namely $X_\ep$ is uniformly equicontinuous in the compact set $[0,T]\times K$. Therefore, by virtue of Ascoli's theorem this combined with \eqref{K'} and \eqref{cae} implies that up to a subsequence $X_\ep$ converges uniformly in $[0,T]\times K$ to a solution $X$ to
\[
\forall\,t\in[0,T],\ \forall\,x\in K,\quad X(t,x)=x+\int_0^t a(X(s,x))\,ds,
\]
{\em i.e.} $X$ is the flow associated with the vector field $a$. Since $a$ belongs to $C^1_b(\RR^N)$, the flow $X$ is uniquely determined (see, {\em e.g.}, \cite[Sec.\,17.4]{HSD}).
Therefore, the whole sequence $X_\ep$ converges uniformly to $X$ in $[0,T]\times K$.
Moreover, by the differentiability of the flow (see, {\em e.g.}, \cite[Sec.\,17.6]{HSD}) we have
\beq\label{DX}
\forall\,t\in[0,T],\ \forall\,x\in K,\quad D_xX_\ep(t,x)=I_N+\int_0^t D_x X_\ep(s,x)\,D_x a_\ep(X_\ep(s,x))\,ds,
\eeq
which using \eqref{cae}, \eqref{K'} and Gronwall's inequality implies that there exists a constant $c>0$ such that
\[
\forall\,t\in[0,T],\ \forall\,x\in K,\quad |D_x X_\ep(t,x)|\leq |I_N|\,e^{ct}.
\]
Therefore, convergences \eqref{cWe} hold true.
\par
On the other hand, by the Liouville formula associated with equation \eqref{DX} and estimate \eqref{aeA} we get that there exists a constant $c>1$ such that
\beq\label{edetXe}
\forall\,t\in[0,T],\ \forall\,x\in K,\quad c^{-1}\leq \det\left(D_xX_\ep(t,x)\right)=\exp\left(\int_0^t (\div\,a_\ep)(X_\ep(s,x))\,ds\right)\leq c,
\eeq
which implies the existence of a constant $C>0$ such that for any $t\in[0,T]$ and $x\in K$,
\[
\ba{l}
\dis \big|\det\,(D_xX_\ep(t,x))-\det\,(D_xX(t,x))\big|
\\ \ecart
\dis \leq C\int_0^T\!\!|\div\,a_\ep-\div\,a|(X_\ep(s,x))\,ds+C\int_0^T\!\!\big|(\div\,a)(X_\ep(s,x))-(\div\,a)(X(s,x))\big|\,ds.
\ea
\]
Hence, using successively Jensen's inequality with respect to the integral in $s$, Fubini's theorem and the change of variables $y=X_\ep(s,x)$ together with \eqref{K'} and \eqref{edetXe}, it follows that there exists a constant $C>0$ such that for any $t\in[0,T]$,
\[
\ba{l}
\dis \big\|\det\left(D_x X_\ep(t,\cdot)\right)-\det\left(D_x X(t,\cdot)\right)\big\|_{L^q(K)}
\\ \ecart
\leq \dis C\,\|\div\,a_\ep-\div\,a\|_{L^q(K')}+C\,\sup_{[0,T]\times K}\big|(\div\,a)(X_\ep)-(\div\,a)(X)\big|.
\ea
\]
This combined with convergence \eqref{divae} and the uniform convergence of $X_\ep$ to $X$ in the compact set $[0,T]\times K$ implies the convergence \eqref{cthe} of $\theta_\ep=\det (D_xX_\ep(t,\cdot))$.
\par
Finally, let $u_\ep^0$ be a sequence in $L^p(\RR^N)$ converging weakly to $u^0$ in $L^p(\RR^N)$. By virtue of Theorem~\ref{thm.hom} combined with Remark~\ref{rem.vsiu} and recalling that $\sigma_\ep=1$, the sequence $u_\ep$ converges weakly in $L^p(\RR^N)$ to a solution $u$ to the equation \eqref{equ}  where $\sigma_0=1$ and by \eqref{xi0cofDW}
\[
\xi_0={\rm Cof}\,(D_x X(t,\cdot))\,e_1\quad\mbox{in }\RR^N.
\]
Proposition~\ref{pro.ae} is thus proved. \cqfd

\end{document}